\documentclass{IEEEtran}
\usepackage{cite}
\usepackage{amsmath,amssymb,amsfonts}
\usepackage{algorithmic}
\usepackage{graphicx}
\usepackage{algorithm,algorithmic}
\usepackage{hyperref}
\usepackage{textcomp}
\usepackage{mathrsfs}
\usepackage{bbm}
\usepackage{amsthm, amsfonts, outlines, mathtools, outlines, cancel, enumitem, float, stackengine, graphicx, romannum, relsize, booktabs, siunitx, dsfont, makecell}

\setenumerate[1]{label=(\roman*)}

\newtheorem{theorem}{Theorem}
\newtheorem{example}{Example}
\newtheorem{assumption}{Assumption}
\newtheorem{lemma}{Lemma}

\newtheorem{propst}{Proposition}
\newtheorem{envdef}{Definition}

\allowdisplaybreaks

\DeclareMathSizes{10}{9}{6}{5}

\newcommand{\abs}[1]{\lvert#1\rvert}%
\newcommand{\babs}[1]{\big\lvert#1\big\rvert}%
\newcommand{\norm}[1]{\lVert#1\rVert}
\newcommand{\bnorm}[1]{\big\lVert#1\big\rVert}

\newcommand{\expt}[2]{\mathbb{E}_{#1}[#2]}
\newcommand{\bexpt}[2]{\mathbb{E}_{#1}\big[#2\big]}
\newcommand{\Bexpt}[2]{\mathbb{E}_{#1}\Big[#2\Big]}

\newcommand\ubar[1]{\stackunder[1.2pt]{$#1$}{\rule{.8ex}{.075ex}}}

\newcommand{\bsb}[1]{\boldsymbol{#1}}
\newcommand{\idty}{\text{Id}}
\newcommand{\probb}{\mathbb{P}}
\newcommand{\mcup}{\,\mathlarger{\mathlarger{\cup}}\,}
\newcommand{\mcap}{\,\mathlarger{\mathlarger{\cap}}\,}

\newcommand{\rset}[2]{\mathbb{R}^{#1}_{#2}}

\newcommand{\nset}[2]{\mathbb{N}^{#1}_{#2}}

\DeclareMathOperator{\argmin}{argmin}

\DeclareMathOperator{\cls}{cl}

\DeclareMathOperator{\ind}{ind}
\DeclareMathOperator{\spans}{Span}

\usepackage{accents}

\DeclareMathOperator{\proj}{P}

\newcommand{\Tblue}[1]{\textcolor{blue}{#1}}

\newif\ifproceeding
\newif\ifarxiv

\proceedingfalse
\arxivtrue

\def\BibTeX{{\rm B\kern-.05em{\sc i\kern-.025em b}\kern-.08em
    T\kern-.1667em\lower.7ex\hbox{E}\kern-.125emX}}
\begin{document}

\title{Offline Learning of Decision Functions in Multiplayer Games with Expectation Constraints}

\author{Yuanhanqing Huang$^{1}$ and Jianghai Hu$^{1}$
\thanks{This work was supported by the National Science Foundation under Grant No. 2014816 and No.2038410. }
\thanks{$^{1}$The authors are with the Elmore Family School of Electrical and Computer Engineering, Purdue University, West Lafayette, IN, 47907, USA 
        {\tt\small \{huan1282, jianghai\}@purdue.edu}}%
}

\maketitle

\begin{abstract}
We explore a class of stochastic multiplayer games where each player in the game aims to optimize its objective under uncertainty and adheres to some expectation constraints. 
The study employs an offline learning paradigm, leveraging a pre-existing dataset containing auxiliary features. 
While prior research in deterministic and stochastic multiplayer games primarily explored vector-valued decisions, this work departs by considering function-valued decisions that incorporate auxiliary features as input. 
We leverage the law of large deviations and degree theory to establish the almost sure convergence of the offline learning solution to the true solution as the number of data samples increases. 
\end{abstract}


\section{Introduction}
Multiplayer (continuous) games model the scenario where a group of interacting rational players choose decisions from continuous decision spaces to optimize their objectives, which are inherently intertwined with the decisions made by others. 
Typical applications include thermal control in smart buildings \cite{jiang2021game, jiang2023game}, path planning for self-driving cars \cite{wang2021game}, supply chain management \cite{hall2024game}, etc. 
In multiplayer games, Nash equilibrium (NE) serves as a pivotal solution concept, characterizing a stable decision profile where no rational player can benefit via unilaterally deviating from it \cite{nash1950equilibrium}. 

Existing work, such as \cite{tatarenko2023accelerating, belgioioso2022distributed, bianchi2022fast}, proposes algorithms for finding NEs and their variants in deterministic multiplayer games, assuming players possess full/partial knowledge of games or access to precise oracles for partial-gradient/first-order information. 
In situations where first-order information is unavailable, recent studies \cite{huang2023zeroth, huang2023global, tatarenko2024payoff, meng2023decentralized} are dedicated to solving NE problems with only bandit oracles/payoff information. 
Real-world setups often involve objective functions affected by uncertainty or random variables, prompting players to optimize expectations of these functions. 
The prevalence of uncertainty necessitates investigation of stochastic NE problems, which has been explored in \cite{lei2022stochastic, huang2022distributed, huang2022distlearn}. 
Most of the previously discussed methods seek to locate the solution decision vector within $\rset{n}{}$. 
Even in stochastic settings, the solution remains in its vectorial form, ensuring that no unilateral deviation can yield improved expected objective values. 
Nevertheless, in the current data-driven landscape, decision-makers can observe auxiliary features before making decisions, despite the unobservability of uncertain quantities at the time of decision-making. 
Inspired by optimizer prediction methods discussed in \cite{bertsimas2020predictive, bertsimas2022data}, we examine the scenario where each player employs a decision function to obtain the optimal decision vector by incorporating realized auxiliary features as inputs. 

Our exploration centers on offline learning, where players leverage a static, finite dataset to formulate their decision functions, as opposed to the continuous and ongoing data stream in the online setting \cite{mertikopoulos2023unified, huang2023zeroth, huang2023global, tatarenko2024payoff, hsieh2021adaptive, meng2023decentralized}. 
Offline learning serves as the preferred paradigm in safety-critical domains demanding rigorous data quality and model performance assurance such as collision avoidance of self-driving cars, medical diagnosis systems, and automated fraud account enforcement. 

\textit{Contributions: }
We examine the offline learning paradigm for multiplayer continuous games under uncertainty. 
While most existing studies are centered around deriving a static equilibrium, our objective is to acquire a set of decision functions in finite-dimensional function spaces that yield varying decision profiles/variables depending on observed auxiliary features, while preserving a sense of stability analogous to that of the NE. 
Additionally, the individual optimization problem for each player incorporates nonlinear expectation constraints, enabling the modeling of scenarios such as obstacle avoidance and risk-related measures, e.g., conditional value-at-risk (CVaR)\cite{rockafellar2002conditional}. 
In contrast to the online learning paradigm, where convergence is often established through Doob's martingale convergence theorems and their variants, the convergence in offline learning with expectation constraints follows a different proof path by comparing the solution sets of two constrained optimization problems. 
Our analysis establishes that finite-sample versions of pseudo-gradients and feasible sets closely approximate their expected counterparts within specified error bounds with high probability. 
In particular, the probability of deviations decreases exponentially as the size of the data sample grows. 
Furthermore, we demonstrate that solutions to the variational inequalities (VIs) constructed using the approximated pseudo-gradient and feasible sets converge almost surely (a.s.) to the true solutions. 
\Tblue{The validity of our method is demonstrated via a multi-account portfolio optimization problem \cite[Appendix~C]{huang2024offline}}. 
To the best of our knowledge, this work initiates the foray into applying optimizer prediction methods \cite{bertsimas2022data} to multiplayer games where decisions are represented as functions. 

\textit{Related Work: }
In their work \cite{peters2023online}, Peters et al. explored the efficacy of the offline learning paradigm to find open-loop NEs in games characterized by smooth open-loop dynamics through comprehensive empirical studies without convergence proof. 
The sample average approximation (SAA) method for single agent stochastic programming (SP) problems with expectation constraints proposed by Wang and Ahmed \cite{wang2008sample} exhibits guaranteed convergence of the feasible set approximation. 
However, the convergence properties of candidate solutions remain unexplored. 
Lan and Zhou \cite{lan2020algorithms} explored two single-agent SP problems with expectation constraints.
Their work introduced a novel stochastic approximation (SA) type algorithm, operating in an online learning manner. 
Koppel et al. \cite{koppel2019projected} investigated single-agent SP problems with expectation constraints, where decisions are functions in a reproducing kernel Hilbert space (RKHS). 
Their online learning algorithm is based on a projected stochastic primal-dual scheme and exhibits a critical no-regret property. 
Iusem et al. \cite{iusem2017extragradient} proposed an online extragradient method for stochastic VIs with sample average approximation (SAA), which is applicable to solving multiplayer stochastic games subject to expectation constraints. 
In a seminal work \cite{bertsimas2022data}, Bertsimas and Koduri proposed two regression-based methods, leveraging RKHS, to tackle constraint-free SP problems. 
By incorporating a set of historical data with auxiliary information, their offline learning approaches ensure convergence of cost values to optimal values in probability as data size grows.  

We remark that the results of this work can be extended to multivariate decision functions and multi-dimensional expectation constraint functions, while we restrict ourselves to scalar decisions and constraint functions for presentation simplicity. 
\ifproceeding
The complete proof of this paper can be found in \cite{huang2024offline} due to space constraints.
\fi

\textit{Notations}: For an integer $n \in \nset{}{+}$, denote $[n] \coloneqq \{1, \ldots, n\}$. 
For $x, y \in \rset{}{}$, $x \wedge y$ represents the minimum value between them. 
Given sets $S_1, S_2, \ldots, S_n$, let $\prod_{i \in [n]} S_i$ denote their Cartesian product. 
Given vectors $v_1, v_2, \ldots, v_n$, let $[v_i]_{i \in [n]}$ and $[v_1; v_2; \cdots; v_n]$ denote their vertical stack; let $[v_1, v_2, \ldots, v_n]$ denote their horizontal stack. 
For a vector $v$, let $[v]_i$ denotes its $i$-th element. 
We let $\norm{v}_2 = \sqrt{x_1^2 + \cdots + x_m^2}$ denote the $\ell_2$ norm of the vector $v = [x_1, \ldots, x_m]^T$.  
For a set $\mathcal{N}$, we use $\overline{\mathcal{N}}$ and $\cls \mathcal{N}$ to denote its closure, while using $\partial \mathcal{N}$ to denote its boundary.  
For a closed and convex set $\mathcal{S} \subseteq \rset{m}{}$ and $x_* \in \rset{m}{}$, $P_{\mathcal{S}}(x_*) = \argmin_{x \in \mathcal{S}} \norm{x - x_*}_2$. Additional notations are provided in Table~\ref{tab:symbol} to provide clarity and readability.

\begin{table}
    \centering
    \begin{tabular}{cc}
        \toprule
        Symbol & Meaning  \\
        \midrule
        $\hat{z}^i$ & Decision function of player $i$  \\
        $z^i$ & Output value of decision function $\hat{z}^i$  \\
        $a^i$ & Coefficients to determine $\hat{z}^i$ \\
        $J^i$, $h^i$ & \makecell{Objective function and constraint function of player $i$ \\ with the realized decision function value $z^i$ as input } \\
        $\mathbf{J}^i$, $\mathbf{h}^i$ & \makecell{Objective function and constraint function of player $i$ \\ with the coefficient $a^i$ as input } \\
        $F$ & Pseudo-gradient of the game with $z^i$'s of all players as input \\
        $\mathbf{F}$ & Pseudo-gradient of the game with $a^i$'s of all players as input \\
        $\mathbb{F}$ & Expected-value version of $\mathbf{F}$ \\
        $\mathbf{F}_{\omega, n}$ & The finite-sample approximation of $\mathbb{F}$ \\
        \bottomrule
    \end{tabular}
    \caption{Key Notations Used in the Paper}
    \label{tab:symbol}
    \vspace*{-\baselineskip}
\end{table}

\section{Problem Setup and Preliminaries}

\subsection{Problem Setup and Finite-Sample Approximations}
We consider a game $\mathcal{G}$ consisting of $N$ rational and self-interested players indexed by $[N]$. 
Unlike the common setup where each player aims to find a decision variable $z^i \in \rset{}{}$ that optimizes its local objective $J^i$, in this work, each player $i \in [N]$ aims to choose a decision function $\hat{z}^i$ that optimizes its local objective function $J^i$ in its expected value. 
The decision function $\hat{z}^i: \mathcal{X} \to \rset{}{}$ determines the decision variable $z^i = \hat{z}^i(X)$ of player $i$ based on some auxiliary and \textit{observable} random factors, captured by the $\mathcal{X}$-valued r.v. $X:\Omega \to \mathcal{X}$, where $(\Omega, \mathcal{F}, \mathbb{P})$ is the underlying probability space. 
Once the decision variables of all players are determined, the local objective of player $i$ is given by $J^i(z^i;y, z^{-i})$.  
Here, $y \in \mathcal{Y} \subseteq \rset{n_y}{}$ is the outcome of a r.v. $Y: \Omega \to \mathcal{Y}$ capturing all the {\it unobservable} random factors affecting the local objectives; $z^{-i} \in \rset{N-1}{}$ are the decision variables made by all players except $i$, which can be thought of as the realized values of $\hat{z}^{-i}: \mathcal{X} \to \rset{N-1}{}$, the stack of the decision functions of these players.  
Note that in $J^i(z^i; y, z^{-i})$, $z^i$ is the function's input, while $y$ and $z^{-i}$ are treated as parameters. 

\begin{assumption}\label{asp:set}
The sets $\mathcal{X} \subseteq \rset{n_x}{}$ and $\mathcal{Y} \subseteq \rset{n_y}{}$ are compact. 
For each $i \in [N]$, given any $y \in \mathcal{Y}$ and $z^{-i}$, the objective function $J^i(\cdot; y, z^{-i})$ is continuously differentiable, whose derivative regarding player $i$'s decision variable is denoted as $\partial J^i(\cdot; y, z^{-i})$.  
\end{assumption}
\begin{example}
(Optimal Production Plans of Vendors) 
Consider a scenario involving $N$ manufacturers, where each manufacturer $i \in [N]$ is required to make decisions on the production quantity of commodities before the disclosure of the market demand $Y = [Y^i]_{i \in [N]}$. 
During decision-making, manufacturers can access factors denoted as $X$ that influence demand, which encompasses elements such as seasonal variations, consumer preferences, economic conditions, etc. 
Following this, each manufacturer $i$ devises its decision function $\hat{z}^i$, mapping the observable features $X$ to the corresponding production quantity $z^i$. 
Its objective function is given by 
\begin{align*}
J^i(z^i; y, z^{-i}) = \big(d^i + \alpha \sum_{j \in [N]} z^j\big)z^i - r^i (z^i \wedge y^i)
\end{align*}
for the decision variables $z^i \in \rset{}{+}$ and $z^{-i} \in \rset{N-1}{+}$. 
Here, 
$\alpha$ describes the influence of the collective production plan on the prices of raw materials; 
the unit base cost of production and the unit revenue are denoted as $d^i$ and $r^i$. 
\end{example}
We restrict our effort in this work to the case of finite-dimensional function spaces, exemplified by linear spaces of homogeneous polynomials and polynomials with degrees not exceeding a specified maximum \cite[Sec.~2.1]{cucker2007learning}. 
Given a finite set of arbitrary continuous functions $\Phi^d \coloneqq \{ \phi_1, \ldots, \phi_d \} \subseteq \mathcal{C}(\mathcal{X})$, we consider the function space $\mathcal{H}$ formed by their linear span. 
For any $\hat{z}^i = \sum_{l \in [d]} a^i_l\phi_l$ and $\hat{w}^i = \sum_{l \in [d]} b^i_l\phi_l$ ($a^i, b^i \in \rset{d}{}$), define the inner product as $\langle \hat{z}^i \mid \hat{w}^i \rangle_{\mathcal{H}} = \langle [a^i_l]_{l \in [d]}, [b^i_l]_{l \in [d]}\rangle = \sum_{l \in [d]}[a^i_l]_{l \in [d]} [b^i_l]_{l \in [d]}$. 
For coefficient vector $a^i \coloneqq [a^i_l]_{l \in [d]}$, we assume $a^i \in \cls \mathbb{B}_R$, where $\mathbb{B}_R \coloneqq \{a^i \in \rset{d}{} \mid \norm{a^i}_2 < R\}$. 
It follows that the image set $\hat{\mathcal{Z}}^i \coloneqq \mcup \{\hat{z}^i(\mathcal{X}) \mid a^i \in \cls \mathbb{B}_R\}$ of each decision function $\hat{z}^i$ is a bounded set, due to the continuity of each $\phi_l$, the compactness of $\mathcal{X}$, and the boundedness of each $a^i$. 

In addition to the bounded condition posited above, we assume that the decision function $\hat{z}^i$ should satisfy an expectation constraint $\expt{}{h^i(\hat{z}^i(X); Y)} \leq 0$, where $h^i: \hat{\mathcal{Z}}^i \times \mathcal{Y} \to \rset{}{}$. 
Motivating examples include risk-aware learning via CVaR \cite{rockafellar2002conditional} and semi-supervised learning \cite[Ch.~5]{chapelle2006semi}\cite{bennett1998semi}, where the objective function is used to measure the fidelity to labeled data, while constraints ensure that model maintains proximity for unlabeled data with similar features, \Tblue{etc \cite{lan2020algorithms}.} 
Then the optimization regarding the function $\hat{z}^i$ can be reformulated as a problem involving the coefficients $a^i$'s and presented as follows:  
\begin{align}\label{eq:prob-a}
\begin{split}
& a^i_\star = \argmin_{a^i \in \mathbb{A}^i} \bexpt{X, Y}{\mathbf{J}^i\big(a^i; X, Y, a^{-i}\big)} + \frac{\lambda}{2}\norm{a^i}^2_2 \\
& \mathbb{A}^i \coloneqq \{a^i \in \cls\mathbb{B}_R \mid \mathds{h}^i(a^i) = 
 \expt{X, Y}{\mathbf{h}^i(a^i; X, Y)} \leq 0\},  
\end{split}
\end{align}
where $\mathbf{J}^i$ is a reformulated expression of $J^i$, taking $a^i$ as input, i.e., $\mathbf{J}^i(a^i; x, y, a^{-i}) = J^i(\hat{z}^i(x); y, \hat{z}^{-i}(x))$, with $\hat{z}^i(x) = \sum_{l \in [d]}a^i_l\phi_l(x)$ and $\hat{z}^{-i}(x) = [\sum_{l \in [d]}a^j_l\phi_l(x)]_{j \in [N]\backslash\{i\}}$; similarly for $\mathbf{h}^i$ and $h^i$, which take $a^i$ and $\hat{z}^i(x)$ as input, respectively. 
The term $\frac{\lambda}{2}\norm{a^i}^2_2$ serves as a regularization term that penalizes overly complex $\hat{z}^i$ to prevent overfitting. 
Furthermore, since the joint distribution of the r.v.'s $X$ and $Y$ remains unknown to us, instead of addressing the expected-value problems in \eqref{eq:prob-a} directly, we focus on their scenario-based or empirical counterparts. 
Let $\{(X^k, Y^k)\}_{k \in [n]}$ be an i.i.d. sequence of $n$ r.v. pairs defined in the probability space $(\Omega, \mathcal{F}, \mathbb{P})$, each with the same distribution as $(X, Y)$. 
For any sample path $\omega \in \Omega$, the dataset $\mathcal{S}_{\omega, n}$ consists of the realized values of the sequence $\{(X^k_\omega, Y^k_\omega)\}_{k \in [n]}$. 
A finite-sample approximation of \eqref{eq:prob-a} can be formulated as: 
\begin{align}\label{eq:prob-fs}
\begin{split}
& a^i_{\omega,n} = \underset{a^i \in A^i_{\omega,n}}{\argmin} \frac{1}{n}\sum_{k \in [n]}\mathbf{J}^i\big(a^i; X^k_{\omega}, Y^{k}_{\omega}, a^{-i}\big) + \frac{\lambda}{2}\norm{a^i}^2_2 \\
& A^i_{\omega,n} \coloneqq \{a^i \in \cls\mathbb{B}_R \mid \mathbf{h}^i_{\omega, n}(a^i) = 
\frac{1}{n}\sum_{k \in [n]} \mathbf{h}^i(a^i; X^k_{\omega}, Y^{k}_{\omega}) \leq 0\}. 
\end{split}
\end{align}

\subsection{Standing Assumptions and Solution Concept}
We proceed to introduce the basic version of pseudo-gradient operator, three variants of it, and three additional assumptions essential for our later solution analysis. 
Denote the feasible set of the stacked decisions as $\hat{\mathcal{Z}} \coloneqq \prod_{i \in [N]} \hat{\mathcal{Z}}^i$. 
The pseudo-gradient operator $F: \hat{\mathcal{Z}} \to \rset{N}{}$ is defined as the direct product of each individual derivative: 
\begin{alignat}{1}
\begin{split}
F(z; y) &= \big[ \partial J^i(z^i; y, z^{-i}) \big]_{i \in [N]}, \forall y \in \mathcal{Y}. 
\end{split}
\end{alignat}
Under Assumption~\ref{asp:set} and within the space $\mathcal{H}$, the gradient of $\mathbf{J}^i$ w.r.t. the undetermined coefficients $a^i_l$ inside $\hat{z}^i = \sum_{l \in [d]} a^i_l\phi_l$ is $\nabla \mathbf{J}^i(a^i; x, y, a^{-i}) = \partial J^i(\hat{z}^i(x); y, \hat{z}^{-i}(x)) \cdot \Phi^d(x)$ where 
$\Phi^d(x) \coloneqq [\phi_l(x)]_{l \in [d]} \in \rset{d}{}$ by chain rule. 
Then we can define the scenario-based pseudo-gradient operators $\mathbf{F}: \rset{Nd}{} \times \mathcal{X} \times \mathcal{Y} \to \rset{Nd}{}$ in terms of the coefficients $\bsb{a} \coloneqq [a^i]_{i \in [N]}$ as 
\begin{alignat}{1}
\begin{split}
\mathbf{F}(\bsb{a}; x, y) &= \big[ \nabla \mathbf{J}^i(a^i; x, y, a^{-i}) \big]_{i \in [N]}\\
& = \big[ \partial J^i(\hat{z}^i(x); y, \hat{z}^{-i}(x)) \cdot \Phi^d(x)\big]_{i \in [N]}. 
\end{split}
\end{alignat}
Upon defining the scenario-based operators $F$ and $\mathbf{F}$, we introduce its expected-value version $\mathbb{F}: \rset{Nd}{} \to \rset{Nd}{}$ by replacing the realized values $x$ and $y$ with their random counterparts $X$ and $Y$
\begin{align}
\begin{split}
& \mathbb{F}(\bsb{a}) = \mathbb{F}^0(\bsb{a}) + \lambda\bsb{a}, \mathbb{F}^0(\bsb{a}) = \bexpt{X, Y}{\mathbf{F}(\bsb{a}; X, Y)}. 
\end{split}
\end{align}
where the superscript $^0$ is used to denote the case where there is no regularization, i.e., $\lambda = 0$. 
The interchangeability of derivative and expectation can be justified by Assumption~\ref{asp:set}, the mean-valued theorem, and the dominated convergence theorem. 
As discussed previously, the exact formula of $\mathbb{F}$ is inaccessible. 
In this work, we will approach it via its finite-sample empirical approximation $\mathbf{F}_{\omega, n}: \rset{Nd}{} \to \rset{Nd}{}$: 
\begin{align}
\begin{split}
& \mathbf{F}_{\omega, n}(\bsb{a}) = \mathbf{F}^0_{\omega, n}(\bsb{a}) + \lambda\bsb{a}, \mathbf{F}^0_{\omega, n}(\bsb{a}) = \frac{1}{n}\sum_{k \in [n]}\mathbf{F}(\bsb{a}; X^k_{\omega}, Y^k_{\omega}). 
\end{split}
\end{align}
Denote the stack of constraint functions as $h \coloneqq [h^i]_{i \in [N]}$; similarly for $\mathbf{h}$, $\mathds{h}$, and $h_{\omega,n}$. 
While Assumptions~\ref{asp:set} stipulates the continuous differentiability of each objective function, Assumption~\ref{asp:cont} introduces Lipschitz continuity for the pseudo-gradient operator and constraint function, which are common setup in work related to game theory and NE seeking methods. 

\begin{assumption}\label{asp:cont}
For any $y \in \mathcal{Y}$, the maps $F$ and $h^i$ ($\forall i \in [N]$) are Lipschitz continuous with constants $L_{F}$ and $L_{h^i}$ on $\hat{\mathcal{Z}}$ and $\hat{\mathcal{Z}}^i$, respectively. 
In other words, for any $y \in \mathcal{Y}$, $z_1 \coloneqq [z^i_1]_{i \in [N]} \in \hat{\mathcal{Z}}$, and $z_2 \coloneqq [z^i_2]_{i \in [N]}  \in \hat{\mathcal{Z}}$, we have
\begin{align}
\begin{split}
& \norm{F(z_1; y) - F(z_2; y)}_2 \leq L_{F} \norm{z_1 - z_2}_2, \\
& \abs{h^i(z^i_1; y) - h^i(z^i_2; y)} \leq L_{h} \abs{z^i_1 - z^i_2}, \forall i \in [N]. 
\end{split}
\end{align}
Moreover, there exists $z \in \hat{\mathcal{Z}}$ such that the images $F(z; \mathcal{Y}) \subseteq \rset{N}{}$ and $h^i(z^i; \mathcal{Y}) \subseteq \rset{}{}$ ($\forall i \in [N]$) are bounded. 
\end{assumption}
Since each $\phi_l$ is continuous, we can derive the more explicit Lipschitz conditions $\norm{\mathbf{F}(\bsb{a}_1; x, y) - \mathbf{F}(\bsb{a}_2; x, y)}_2 \leq L_{\mathbf{F}} \norm{\bsb{a}_1 - \bsb{a}_2}_2$ and $\abs{\mathbf{h}^i(a^i_1; x, y) - \mathbf{h}^i(a^i_2; x, y)}_2 \leq L_{\mathbf{h}}\norm{a^i_1 - a^i_2}_2$ ($\forall i \in [N]$), for coefficients $\bsb{a}_1 = [a^i_1]_{i \in [N]}$ and $\bsb{a}_2 = [a^i_2]_{i \in [N]}$ associated with $z_1$ and $z_2$. 
Moreover, since $\hat{\mathcal{Z}} \subseteq \rset{N}{}$ is bounded, by the continuity posited above, there exist constants $M_F$ and $M_h$ such that $\sup_{z \in \hat{\mathcal{Z}}, y \in \mathcal{Y}} \norm{F(z; y)}_\infty < M_F$ and $\sup_{z \in \hat{\mathcal{Z}}, y \in \mathcal{Y}} \norm{h(z; y)}_\infty < M_h$. 

\begin{assumption}\label{asp:mono}
For any $y \in \mathcal{Y}$, $F(\cdot; y): \hat{\mathcal{Z}} \to \rset{N}{}$ is monotone, that is, $\forall z_1, z_2 \in \hat{\mathcal{Z}}$, $\langle F(z_1; y) - F(z_2; y), z_1 - z_2\rangle \geq 0$. 
For each $i \in [N]$ and $\forall y \in \mathcal{Y}$, $h^i(\cdot; y): \hat{\mathcal{Z}}^i \to \rset{}{}$ is convex. 
\end{assumption}

By definition, we have $\langle \mathbf{F}(\bsb{a}_1; x, y) - \mathbf{F}(\bsb{a}_2; x, y), \bsb{a}_1 - \bsb{a}_2\rangle 
= \langle F(\hat{z}_1(x); y) \otimes \Phi^d(x) - F(\hat{z}_2(x); y) \otimes \Phi^d(x), \bsb{a}_1 - \bsb{a}_2\rangle
= \langle F(\hat{z}_1(x); y) - F(\hat{z}_2(x); y), \hat{z}_1(x) - \hat{z}_2(x) \rangle$, where $\otimes$ denotes the Kronecker product on vectors.  
It can be derived that $\mathbf{F}$ and $\mathbb{F}$ inherit the monotonicity of $F$. 
In the same vein, $\mathbf{h}^i$ and $\mathds{h}^i$ inherit the convexity of $h^i$ for each $i \in [N]$. 
To ensure non-emptiness of the approximate feasible set $A^i_{\omega, n}$ with high probability when there are adequate samples, we posit the following assumption in the expected-value sense, which is analogous to Slater's condition.
\begin{assumption}\label{asp:feasib}
There exists an $\varepsilon_h > 0$ \Tblue{s.t.} for every $i \in [N]$, $\{a^i \in \cls \mathbb{B}_R \mid \mathds{h}^i(a^i) \leq -\varepsilon_h\} \neq \varnothing$. 
\end{assumption}

Regarding the solution concept, our focus centers on (function) critical points (CPs), whose vector counterparts have been introduced in the previous literature \cite{huang2023zeroth, huang2023global}. 

\begin{envdef}
In the function space $\mathcal{H} = \spans\Phi^d$, a collective decision function $\hat{z}_\star: \mathcal{X} \to \rset{N}{}$ given by $\hat{z}_\star = [\hat{z}^i_\star]_{i \in [N]}$ and $\hat{z}^i_\star = \sum_{l \in [d]} a^i_{l\star} \phi_l$ is a critical point of the game $\mathcal{G}$ if the vector $\bsb{a}_\star = \big[[a^i_{l\star}]_{l \in [d]}\big]_{i \in [N]}$ is a solution to the following variational inequality $\text{VI}(\mathbb{A}, \mathbb{F})$: 
\begin{align}
\langle \mathbb{F}(\bsb{a}_\star), \bsb{a} - \bsb{a}_\star \rangle \geq 0, \forall \bsb{a} \in \mathbb{A} \coloneqq \prod\nolimits_{i \in [N]} \mathbb{A}^i. 
\end{align}
\end{envdef}
When every $J^i(\cdot; y, z^{-i})$ is convex for any feasible $y$ and $z^{-i}$, the solution concept CPs align with NEs. 
It implies that no rational player can enhance its expected objective value by unilaterally altering the coefficients of its decision function. 

\section{Large Deviation Results for Finite-Sample Approximations} 

\subsection{The Theory of Large Deviations} 

The theory of large deviations (LD) \cite[Sec.~2.2]{kleywegt2002sample}\cite{dembo2009large} is pivotal in our proof, demonstrating that the probability of inaccurate approximation decays exponentially fast to zero. 
Given an $\rset{}{}$-r.v. $X$ with mean zero and variance $\sigma^2$, its LD rate function $I$ is defined as $I(u) = \sup_{t \in \rset{}{}}\{tu - \log M(t)\}$, where $M(t) \coloneqq \expt{}{e^{tX}}$ denotes the moment-generating function of $X$. 
Let $X_1, \ldots, X_K$ be an i.i.d. sequence of replications of $X$ and $\bar{X}_K \coloneqq (1/K)\sum_{k \in [K]} X_k$ their average. 
Then for any $\varepsilon > 0$, the following LD inequality holds
\begin{align}\label{eq:ld}
\probb \{ \abs{\bar{X}_K} > \varepsilon\} \leq 2\exp\big(-KI(\varepsilon)\big). 
\end{align}
Moreover, if $M(t)$ is finite in a neighborhood of $t = 0$, $X$ has finite moments, and Taylor's expansion can be leveraged to build the following lower bound for $I$: 
\begin{align}\label{eq:ld-rate}
\exists \varepsilon_\diamond \in (0, \varepsilon_h) \text{ s.t. } I(\varepsilon) \geq \varepsilon^2 / \big(3 \sigma^2 \big), \forall \varepsilon \leq \varepsilon_\diamond. 
\end{align}
Appendix A contains the proofs relevant to the following subsections.

\subsection{Pseudo-Gradient Approximation}

In this subsection, we start by proving that on a finite set $A$ for coefficients $\bsb{a}$, as the number of data samples $n$ grows, the errors between the finite-sample pseudo-gradient $\mathbf{F}_{\omega, n}(\bsb{a})$ and $\mathbb{F}(\bsb{a})$ decay to zero on a sample subset of $\Omega$ with probability approaching $1$. 

\begin{lemma}\label{le:opt-approx-finite}
Suppose Assumptions~\ref{asp:set} and \ref{asp:cont} hold and the set $A \subseteq \cls\mathbb{B}^N_R$ is finite. 
Then for all  $0 < \varepsilon \leq \varepsilon_\diamond$ with $\varepsilon_\diamond$ defined in \eqref{eq:ld-rate}, we have
\begin{align*}
& \probb\big\{ \omega \mid \bnorm{\mathbf{F}_{\omega, n}(\bsb{a}) - \mathbb{F}(\bsb{a})}_2 \geq \varepsilon, \exists \bsb{a} \in A \big\} \\
& \qquad \leq 2Nd\abs{A} \cdot \exp\Big(-\frac{\varepsilon^2n}{3NdM_F^2}\Big). 
\end{align*}
\end{lemma}

Then we extend the conclusion to the infinite set $\cls\mathbb{B}^N_R \coloneqq \cls\big(\prod_{i \in [N]} \mathbb{B}_R\big)$ via the discretization of compact sets and Lipschitz continuity posited. 

\begin{propst}\label{propst:opt-approx-cmpct}
Suppose Assumptions~\ref{asp:set} and \ref{asp:cont} hold. 
Then for all $0 < \varepsilon \leq \varepsilon_\diamond$ with $\varepsilon_\diamond$ defined in \eqref{eq:ld-rate}, we have 
\begin{align*}
& \probb\big\{ \omega \mid \bnorm{\mathbf{F}_{\omega, n}(\bsb{a}) - \mathbb{F}(\bsb{a})}_2 < \varepsilon, \forall \bsb{a} \in \cls\mathbb{B}^N_R \big\} \\
& \qquad \geq 1 - 2Nd\mathfrak{n}_\nu \exp\Big(-\frac{\varepsilon^2 n}{12NdM_F^2}\Big),  
\end{align*}
where $\mathfrak{n}_\nu \coloneqq (R\sqrt{Nd}/\nu)^{Nd}$ and $\nu = \varepsilon/(4L_{\mathbf{F}})$. 
\end{propst}

For the lower bound derived above, although the term $(Nd\mathfrak{n}_\nu)$ grows at the rate of $O(N(\sqrt{N})^N)$, it remains fixed as $N$ denotes the number of player engaged in the game and we consider games of a fixed size in this paper. 
In comparison, the number $n$ of samples can vary and the term $\exp\Big(-\frac{\varepsilon^2 n}{12NdM_F^2}\Big)$ decreases exponentially as $n$ grows. 
When $n$ is reasonably large, the lower bound derived will be positive and asymptotically converge to 1 as $n \to \infty$. 
The same reasoning extends to the results to be presented in Prop.~\ref{propst:set-approx-cmpct}. 

\subsection{Feasible Set Approximation}

In line with the reasoning in the preceding subsection, we will commence by illustrating that, for finite sets, the constraints formulated by empirical constraint functions effectively approximate the original ones when the sample size $n$ is sufficiently large. 
Before proceeding, for each $i \in [N]$ and any $\varepsilon \in \rset{}{}$, we introduce the following approximate set $\mathbb{A}^i_\varepsilon$: 
\begin{align}
\mathbb{A}^i_\varepsilon \coloneqq \{a^i \in \mathbb{B}_R \mid \mathds{h}^i(a^i) = 
 \expt{X, Y}{\mathbf{h}^i(a^i; X, Y)} \leq \varepsilon \}, 
\end{align}
and let $\mathbb{A}_\varepsilon \coloneqq \prod_{i \in [N]} \mathbb{A}^i_\varepsilon$ and similarly for $A_{\omega, n}$, etc.

\begin{lemma}\label{le:set-approx-finite}
Suppose Assumptions~\ref{asp:set}, \ref{asp:cont}, and \ref{asp:feasib} hold and the set $A \subseteq \cls\mathbb{B}^N_R$ is finite. 
Let $\Tilde{A}_\varepsilon \coloneqq A \mcap \mathbb{A}_\varepsilon$ and $\Tilde{A}_{\omega, n} \coloneqq A \mcap A_{\omega, n}$. 
Then for all $0 < \varepsilon \leq \varepsilon_\diamond$ with $\varepsilon_\diamond$ defined in \eqref{eq:ld-rate}, we have
\begin{align*}
\probb\big\{ \omega \mid \Tilde{A}_{-\varepsilon} \subseteq \Tilde{A}_{\omega, n} \subseteq \Tilde{A}_{\varepsilon} \big\}^c \leq 2N\abs{A}\exp\Big(-\frac{\varepsilon^2n}{3 M_h^2}\Big). 
\end{align*}
\end{lemma}

Subsequently, we will extend the LD results to compact sets. 
\begin{propst}\label{propst:set-approx-cmpct}
Suppose Assumptions~\ref{asp:set}, \ref{asp:cont}, and \ref{asp:feasib} hold. 
Then for all $0 < \varepsilon \leq \varepsilon_\diamond$ with $\varepsilon_\diamond$ defined in \eqref{eq:ld-rate}, we have 
\begin{align*}
& \probb\big\{ \omega \mid \mathbb{A}_{-\varepsilon} \subseteq A_{\omega,n} \subseteq \mathbb{A}_{\varepsilon} \big\} \geq 1 - 2N\mathfrak{n}_\nu \exp\Big(-\frac{\varepsilon^2 n}{12M_h^2}\Big),  
\end{align*}
where $\mathfrak{n}_\nu \coloneqq (R\sqrt{Nd}/\nu)^{Nd}$ and $\nu = \varepsilon/(4L_{\mathbf{h}})$. 
\end{propst}

\subsection{Almost Sure Results}
Consider a sequence of error bounds $(\varepsilon_n)_{n \in \nset{}{+}}$ with each $\varepsilon_n = \varepsilon_\diamond n^{-\alpha}$ for some fixed constant $0 < \alpha < 1/2$. 
Additionally, for each $n \in \nset{}{+}$, select $\nu_n = \varepsilon_n/(4L) = \varepsilon_\diamond n^{-\alpha}/(4L)$ with $L \coloneqq L_F \vee L_h$. 
Define the event $\mathscr{E}_n \coloneqq \big\{ \omega \mid \bnorm{\mathbf{F}_{\omega, n}(\bsb{a}) - \mathbb{F}(\bsb{a})}_2 < \varepsilon_n, \forall \bsb{a} \in \cls\mathbb{B}^N_R \big\} \mcap \big\{ \omega \mid \mathbb{A}_{-\varepsilon_n} \subseteq A_{\omega,n} \subseteq \mathbb{A}_{\varepsilon_n} \big\}$. 

\begin{theorem}\label{thm:as-res}
Suppose Assumptions~\ref{asp:set}, \ref{asp:cont}, and \ref{asp:feasib} hold. Then
\begin{align*}
\probb\big\{ \limsup_{n \to \infty} \mathscr{E}_n^c \big\} = \probb\big\{ \mathscr{E}_n^c \text{ i.o. } \big\} = 0. 
\end{align*}
Then there exists a sample subset $\Tilde{\Omega} \subseteq \Omega$ with $\probb\{\Tilde{\Omega}\} = 1$ where $\forall \omega \in \Tilde{\Omega}$, there exists an index $n_{\omega} \in \nset{}{+}$ such that $\omega \in \mathscr{E}_n$ for all $n > n_{\omega}$. 
\end{theorem}

Following Theorem~\ref{thm:as-res}, starting at index $n_{\omega} + 1$, $\bnorm{\mathbf{F}_{\omega, n}(\bsb{a}) - \mathbb{F}(\bsb{a})}_2 < \varepsilon_n, \forall \bsb{a} \in \cls\mathbb{B}^N_R$ and $\mathbb{A}_{-\varepsilon_n} \subseteq A_{\omega,n} \subseteq \mathbb{A}_{\varepsilon_n}$. 
It enables us to subsequently examine the consistency of solutions for the approximating VIs in the almost sure sense.

\section{Asymptotic Solution Consistency Analysis for Finite-Sample Approximations}

\subsection{Topological Degrees and Existence of Solutions}

The focus of this section is to establish that, with increasing values of $n$, the empirical optimal decision coefficients $\bsb{a}_{\omega,n}$ can arbitrarily approach the true solution $\bsb{a}_\star$ closely on a random sample set with probability $1$. 
(Topological) degree theory for operators serves as a key tool in analyzing the existence and properties of solutions to equations of the form $f(x) = y$, and understanding the effects of modifications to $f$ and $y$. 
The definition, existence, and implications of degree are summarized in the following theorem. 
\begin{theorem}\label{thm:deg}
(\cite[Ch.~1]{deimling2010nonlinear}\cite[Sec.~2.1]{facchinei2003finite}) 
Given the set of tuples
$\mathcal{T} \coloneqq \{(f, \mathcal{X}, y): \mathcal{X} \subseteq \rset{m}{} \text{  open and bounded},  f: \cls \mathcal{X} \to \rset{m}{} \text{continuous}, y \in \rset{m}{} \backslash f(\partial \mathcal{X}) \}$, 
there exists a unique function $\deg: \mathcal{T} \to \mathbb{Z}$ that fulfills the following three axioms \\
(d1) $\deg(\idty, \mathcal{X}, y) = 1$ if $y \in \mathcal{X}$; \\
(d2) $\deg(f, \mathcal{X}, y) = \deg(f, \mathcal{X}_1, y) + \deg(f, \mathcal{X}_2, y)$, if $\mathcal{X}_1$, $\mathcal{X}_2$ are disjoint open subsets of $\mathcal{X}$ and $y \notin f\big( \cls \mathcal{X} \backslash (\mathcal{X}_1 \mcup \mathcal{X}_2) \big)$; \\
(d3) $\deg\big(h(\cdot, t), \mathcal{X}, y(t) \big)$ is independent of $t \in [0,1]$ if $h: \cls\mathcal{X} \times [0,1] \to \rset{m}{}$ and $y:[0,1] \to \rset{m}{}$ are continuous and $y(t) \notin h(\partial \mathcal{X}, t)$ for all $t \in [0, 1]$. 

Moreover, $\deg(f, \mathcal{X}, y) \neq 0 \implies f^{-1}(y) \mcap \mathcal{X} \neq \varnothing$.  
When $x_*$ is an isolated solution to the equation $f(x) = y$, then the degree is common for an arbitrary open neighborhood of $x_*$ that does not contain any other solutions, and such a degree is called the index of $f$ at $x_*$ denoted by $\ind(f, x_*)$. 
\end{theorem}

\subsection{Solution Analysis for Finite-Sample Approximations}

Under Assumption~\ref{asp:mono}, as $F$ is a monotone operator, $\mathbb{F}$ and $\mathbf{F}_{\omega, n}$ are two $\lambda$-strongly monotone operators with the additional regularization terms. 
When $A_{\omega,n}$ is nonempty, $\text{VI}(\mathbb{A}, \mathbb{F})$ and $\text{VI}(A_{\omega,n}, F_{\omega,n})$ admit unique solutions $\bsb{a}_\star$ and $\bsb{a}_{\omega, n}$, respectively \cite[Thm.~2.3.3]{facchinei2003finite}. 
Let $\Tilde{\bsb{a}}_\star = \bsb{a}_\star - \mathbb{F}(\bsb{a}_\star)$ and $\Tilde{\bsb{a}}_{\omega, n} = \bsb{a}_{\omega, n} - \mathbf{F}_{\omega, n}(\bsb{a}_{\omega, n})$. 
As demonstrated in \cite[Prop.~1.5.9]{facchinei2003finite}, $\bsb{a}_\star$ is a solution to $\text{VI}(\mathbb{A}, \mathbb{F})$ if and only if $\mathbb{F}^{\text{nor}}(\Tilde{\bsb{a}}_\star) = 0$, where $\mathbb{F}^{\text{nor}} \coloneqq \mathbb{F} \circ \proj_{\mathbb{A}} + \idty - \proj_{\mathbb{A}}$ is the normal map associated with $\text{VI}(\mathbb{A}, \mathbb{F})$. 
Likewise, for the empirical problem $\text{VI}(A_{\omega, n}, \mathbf{F}_{\omega, n})$, the normal map is given by 
$\mathbf{F}^{\text{nor}}_{\omega, n} \coloneqq \mathbf{F}_{\omega, n} \circ \proj_{A_{\omega, n}} + \idty - \proj_{A_{\omega, n}}$. 
Our analysis will revolve around the index of $\mathbb{F}^{\text{nor}}$ at $\Tilde{\bsb{a}}_\star$, with subsequent emphasis on revealing the analogous behavior of $\mathbb{F}^{\text{nor}}$ and $\mathbf{F}^{\text{nor}}_{\omega, n}$ within the proximity of $\Tilde{\bsb{a}}_\star$. 
The key results are encapsulated in Theorem~\ref{thm:sol-consis}, which demonstrates that the solution to the finite-sample approximation \eqref{eq:prob-fs} converges almost surely to the solution of the original expected value problem \eqref{eq:prob-a}. 
This conclusion follows from the supporting lemmas, Lemma~\ref{le:cont-proj} and Lemma~\ref{le:unit-ind}. 
The proofs related to this subsection can be found in Appendix~B. 

\begin{theorem}\label{thm:sol-consis}
Suppose Assumptions~\ref{asp:set} to \ref{asp:feasib} hold. For any $\omega \in \Tilde{\Omega}$ with $\Tilde{\Omega}$ defined in Theorem~\ref{thm:as-res}, consider an arbitrary neighborhood $\mathcal{N}$ of the solution $\bsb{a}_\star$ to $\text{VI}(\mathbb{A}, \mathbb{F})$.
There exists an $N_{\omega, n} \in \nset{}{+}$ such that for any $n \geq N_{\omega, n}$, the unique solution $\bsb{a}_{\omega, n}$ to $\text{VI}(A_{\omega, n}, \mathbf{F}_{\omega, n})$ satisfies $\bsb{a}_{\omega, n} \in \mathcal{N}$. 
Moreover, it implies that $\lim_{n \to \infty} \norm{\bsb{a}_{\omega, n} - \bsb{a}_\star}_2 = 0$. 
\end{theorem}

Given the expectation constraints in \eqref{eq:prob-a} and \eqref{eq:prob-fs}, Lemma~\ref{le:cont-proj} shows that, for a sufficiently large $n$, the projection operator $\proj_{A_{\omega, n}}$ produces results that closely mimic those of $\proj_{\mathbb{A}}$. 
\begin{lemma}\label{le:cont-proj}
Suppose Assumptions~\ref{asp:set}, \ref{asp:cont}, and \ref{asp:feasib} hold. 
For any $\omega \in \Tilde{\Omega}$ with $\Tilde{\Omega}$ defined in Theorem~\ref{thm:as-res}, $\lim_{n \to \infty} A_{\omega, n} = \mathbb{A}$. 
For an arbitrary $\Tilde{\bsb{a}}_\sharp \in \rset{Nd}{}$ and $\bsb{a}_\sharp \coloneqq \proj_{\mathbb{A}}(\Tilde{\bsb{a}}_\sharp)$, consider an arbitrary open neighborhood $\mathcal{N}$ of $\bsb{a}_\sharp$. 
Then there exists $N_{0, \omega} \in \nset{}{}$ and an open neighborhood $\Tilde{\mathcal{N}}$ of $\Tilde{\bsb{a}}_\sharp$ s.t. 
$\proj_{A_{\omega, n}}\big( \cls \Tilde{\mathcal{N}} \big) \subseteq \mathcal{N}, \forall n > N_{0, \omega}$. 
\end{lemma}

To invoke Theorem~\ref{thm:deg}, Lemma~\ref{le:unit-ind} is dedicated to establishing that for the original problem \eqref{eq:prob-a}, the index of the normal map at the point $\Tilde{\bsb{a}}_\star \coloneqq \bsb{a}_\star - \mathbb{F}(\bsb{a}_\star)$ is non-zero. 
\begin{lemma}\label{le:unit-ind}
Suppose Assumptions~\ref{asp:set} to \ref{asp:mono} hold. 
Then for the unique solution $\bsb{a}_\star$ of $\text{VI}(\mathbb{A}, \mathbb{F})$, the index of $\mathbb{F}^{\text{nor}}$ at $\Tilde{\bsb{a}}_\star \coloneqq \bsb{a}_\star - \mathbb{F}(\bsb{a}_\star)$ is $\ind (\mathbb{F}^{\text{nor}}, \Tilde{\bsb{a}}_\star) = 1$. 
\end{lemma}

\section{Conclusion}
This work investigates stochastic multiplayer games, characterized by expected value objectives and expectation constraints. 
Leveraging historical datasets and auxiliary features, the proposed offline learning method produces a set of function-valued decisions that provably converge a.s. to stable decision functions characterized by the associated variational inequality. 
This work is currently restricted to finite-dimensional function spaces. 
In forthcoming endeavors, we aim to enhance the descriptive capability of function spaces and explore the possibility of extending the results to RKHS with universal kernels. 

\bibliographystyle{IEEEtran}
\bibliography{IEEEabrv,references_abbr}

\appendices 
\section*{Appendix}
\addcontentsline{toc}{section}{Appendix}
\renewcommand{\thesubsection}{\Alph{subsection}}

\newtheorem{appdxlemma}{Lemma}
\newtheorem{appdxtheorem}{Theorem}
\newtheorem{appdxpropst}{Proposition}
\numberwithin{appdxtheorem}{subsection}
\numberwithin{appdxlemma}{subsection} 
\numberwithin{appdxpropst}{subsection}
\numberwithin{equation}{subsection} 

\subsection{Proof for LD Results}
\subsubsection{Proof of Lemma~\ref{le:opt-approx-finite}}
For notational ease, we will occationally omit ``$\omega |$'' in expressions of events such as $\{\omega | f(\omega) \leq 0\} \in \mathcal{F}$. 
For each $\bsb{a} \in A$, we have
$\big\{ \bnorm{\mathbf{F}_{\omega, n}(\bsb{a}) - \mathbb{F}(\bsb{a})}_2 \geq \varepsilon \big\} 
\subseteq \big\{ \babs{[\mathbf{F}_{\omega, n}(\bsb{a}) - \mathbb{F}(\bsb{a})]_l} \geq \varepsilon/\sqrt{Nd}, \exists l \in [Nd] \big\}
= \mcup_{l \in [Nd]} \big\{ \babs{[\mathbf{F}_{\omega, n}(\bsb{a}) - \mathbb{F}(\bsb{a})]_l} \geq \varepsilon/\sqrt{Nd} \big\}$. 
It follows that 
\begin{align*}
& \probb\big\{ \bnorm{\mathbf{F}_{\omega, n}(\bsb{a}) - \mathbb{F}(\bsb{a})}_2 \geq \varepsilon, \exists \bsb{a} \in A \big\} \\
& \leq \sum_{\bsb{a} \in A}\sum_{l \in [Nd]} \probb \big\{ \babs{[\mathbf{F}_{\omega, n}(\bsb{a}) - \mathbb{F}(\bsb{a})]_l} \geq \varepsilon/\sqrt{Nd} \big\} \\
& \leq 2Nd\abs{A} \exp\Big(-n \cdot \frac{\varepsilon^2}{3NdM_F^2}\Big), 
\end{align*}
where in the last line we apply \eqref{eq:ld} and note that the variance satisfies $\bexpt{}{\babs{[\mathbf{F}(\bsb{a}; X, Y)]_l}^2} - [\mathbb{F}(\bsb{a})]_l^2 \leq M_F^2$.

\subsubsection{Proof of Proposition~\ref{propst:opt-approx-cmpct}}
Since $\cls \mathbb{B}^N_R$ is compact, for any $\nu > 0$, we can cover the set $\cls\mathbb{B}^N_R$ with finitely many open balls with radius $\nu$. 
By noting that for any $\bsb{a}_1, \bsb{a}_2 \in \rset{Nd}{}$, $\norm{\bsb{a}_1 - \bsb{a}_2}_\infty \leq \nu/\sqrt{Nd} \implies \norm{\bsb{a}_1 - \bsb{a}_2}_2 \leq \nu$, a suboptimal estimate of the count is $\mathfrak{n}_\nu = (R\sqrt{Nd}/\nu)^{Nd}$. 
Denote the centers of this open cover by $A_\nu \coloneqq \{\bsb{a}^\nu_1, \ldots, \bsb{a}^\nu_{\mathfrak{n}_\nu}\}$. 
Then for an arbitrary $\bsb{a} \in \cls \mathbb{B}^N_R$, there exists an $\bsb{a}^\nu \in A_\nu$ such that $\norm{\bsb{a} - \bsb{a}^\nu}_2 < \nu$. 
Following that, we consider the following decomposition of the error: 
\begin{align*}
\mathbf{F}_{\omega, n}(\bsb{a}) - \mathbb{F}(\bsb{a}) =& F^0_{\omega, n}(\bsb{a}) - F^0_{\omega, n}(\bsb{a}_\nu) + F^0_{\omega, n}(\bsb{a}_\nu) \\
& - \mathbb{F}^0(\bsb{a}_\nu) + \mathbb{F}^0(\bsb{a}_\nu) - \mathbb{F}^0(\bsb{a}). 
\end{align*}
By the Lipschitz continuity of $F$, it is evident that 
\begin{align*}
& \bnorm{F^0_{\omega, n}(\bsb{a}) - F^0_{\omega, n}(\bsb{a}_\nu)}_2 \\
& \leq \frac{1}{n} \sum\nolimits_{k \in [n]} \bnorm{\mathbf{F}(\bsb{a}; X^k_{\omega}, Y^k_{\omega}) - \mathbf{F}(\bsb{a}_\nu; X^k_{\omega}, Y^k_{\omega})}_2 \\
& \leq \frac{1}{n} \sum\nolimits_{k \in [n]} L_{\mathbf{F}}\norm{\bsb{a} - \bsb{a}_\nu}_2  < L_{\mathbf{F}}\nu, \\
& \bnorm{\mathbb{F}^0(\bsb{a}) - \mathbb{F}^0(\bsb{a}_\nu)}_2 \leq \bexpt{}{\bnorm{ \mathbf{F}(\bsb{a}; X, Y) - \mathbf{F}(\bsb{a}_\nu; X, Y) }_2} \\
& \leq \bexpt{}{L_{\mathbf{F}}\norm{\bsb{a} - \bsb{a}_\nu}_2} < L_{\mathbf{F}}\nu. 
\end{align*}
Hence, the following inclusions of events hold: 
\begin{align*}
& \big\{\bnorm{\mathbf{F}_{\omega, n}(\bsb{a}) - \mathbb{F}(\bsb{a})}_2 \geq \varepsilon, \exists \bsb{a} \in \cls\mathbb{B}^N_R\big\} \\
& \subseteq \big\{\bnorm{\mathbf{F}_{\omega, n}(\bsb{a}_\nu) - \mathbb{F}(\bsb{a}_\nu)}_2 + 2L_{\mathbf{F}}\nu \geq \varepsilon, \exists \bsb{a} \in \cls\mathbb{B}^N_R, \\
& \qquad \text{and some } \bsb{a}_\nu \in A_\nu \text{ s.t. } \norm{\bsb{a} - \bsb{a}_\nu} < \nu \big\} \\
& \subseteq \big\{\bnorm{\mathbf{F}_{\omega, n}(\bsb{a}_\nu) - \mathbb{F}(\bsb{a}_\nu)}_2 + 2L_{\mathbf{F}}\nu \geq \varepsilon, \exists \bsb{a}_\nu \in A_\nu \big\} \\
& = \big\{\bnorm{\mathbf{F}_{\omega, n}(\bsb{a}_\nu) - \mathbb{F}(\bsb{a}_\nu)}_2  \geq \varepsilon/2, \exists \bsb{a}_\nu \in A_\nu \big\}, 
\end{align*}
where the last inclusion follows from the condition that $\nu =  \varepsilon/(4L_{\mathbf{F}})$. 
In light of Lemma~\ref{le:opt-approx-finite}, it can be deduced that 
\begin{align*}
& \probb\big\{\omega \mid \bnorm{\mathbf{F}_{\omega, n}(\bsb{a}) - \mathbb{F}(\bsb{a})}_2 \geq \varepsilon, \exists \bsb{a} \in \cls\mathbb{B}^N_R\big\} \\
& \qquad \leq 2Nd\mathfrak{n}_\nu \exp\Big(-\frac{\varepsilon^2 n}{12NdM_F^2}\Big), 
\end{align*}
and the conclusion follows by taking the complement. 

\subsubsection{Proof of Lemma~\ref{le:set-approx-finite}}
We begin by observing that 
\begin{align*}
& \big\{\Tilde{A}_{-\varepsilon} \subseteq \Tilde{A}_{\omega, n} \subseteq \Tilde{A}_{\varepsilon} \big\}^c \\
& = \big\{ \exists \bsb{a} \in \Tilde{A}_{-\varepsilon} \text{ s.t. } \bsb{a} \notin \Tilde{A}_{\omega, n} \big\} \mcup \big\{ \exists \bsb{a} \in \Tilde{A}_{\omega, n} \text{ s.t. } \bsb{a} \notin \Tilde{A}_{\varepsilon} \big\} \\
& \subseteq \big\{ \exists \bsb{a} \in A, i \in [N] \text{ s.t. } \mathds{h}^i(a^i) \leq -\varepsilon, \text{ and } \mathbf{h}^i_{\omega, n}(a^i) > 0\big\} \\
& \qquad \mcup \big\{ \exists \bsb{a} \in A, i \in [N] \text{ s.t. } \mathds{h}^i(a^i) > \varepsilon, \text{ and } \mathbf{h}^i_{\omega, n}(a^i) \leq 0\big\} \\
& \subseteq \big\{ \exists \bsb{a} \in A, i \in [N] \text{ s.t. } \abs{\mathbf{h}^i_{\omega, n}(a^i) - \mathds{h}^i(a^i)} > \varepsilon \big\}. 
\end{align*}
Based on the theory of LD, we can procure the following: 
\begin{align*}
\probb\big\{ \Tilde{A}_{-\varepsilon} \subseteq \Tilde{A}_{\omega, n} \subseteq \Tilde{A}_{\varepsilon} \big\}^c & \leq \sum_{\bsb{a} \in A}\sum_{i \in [N]} \probb\big\{ \abs{\mathbf{h}^i_{\omega, n}(a^i) - \mathds{h}^i(a^i)} > \varepsilon \big\} \\
&  \leq 2 N\abs{A}\exp\Big( - \frac{n\varepsilon^2}{3M_h^2} \Big), 
\end{align*}
where in the last line we apply \eqref{eq:ld} and note that the variance satisfies $\bexpt{}{\babs{\mathbf{h}^i(a^i; X, Y)}^2} - \mathds{h}^i(a^i)^2 \leq M_h^2$. 

\subsubsection{Proof of Proposition~\ref{propst:set-approx-cmpct}}
Employing identical finite covering rationale in Proposition~\ref{propst:opt-approx-cmpct}, we attain a set of centers $A_\nu \coloneqq \{\bsb{a}^\nu_1, \ldots, \bsb{a}^\nu_{\mathfrak{n}_\nu}\}$ with $\mathfrak{n}_\nu = (R\sqrt{Nd}/\nu)^{Nd}$ s.t. for any $\bsb{a} = [a^i]_{i \in [N]} \in \cls \mathbb{B}^N_R$, there exists an $\bsb{a}_\nu = [a^i_{\nu}]_{i \in [N]} \in A_\nu$ satisfying $\norm{\bsb{a} - \bsb{a}_\nu}_2 < \nu$. 
Utilizing the error decomposition arguments yields 
$\mathbf{h}^i_{\omega, n}(a^i) - \mathds{h}^i(a^i) = \mathbf{h}^i_{\omega, n}(a^i) - \mathbf{h}^i_{\omega, n}(a^i_\nu) + \mathbf{h}^i_{\omega, n}(a^i_\nu) - \mathds{h}^i(a^i_\nu) + \mathds{h}^i(a^i_\nu) - \mathds{h}^i(a^i)$. 
Additionally, it follows from the Lipschitz continuity of each $h^i$ that: 
\begin{align*}
& \babs{\mathbf{h}^i_{\omega, n}(a^i) - \mathbf{h}^i_{\omega, n}(a^i_\nu)} \\
& \leq \frac{1}{n}\sum\nolimits_{k \in [n]} \babs{\mathbf{h}^i(a^i; X^k_\omega, Y^k_\omega) - \mathbf{h}^i(a^i_\nu; X^k_\omega, Y^k_\omega)} \\
& \leq \frac{1}{n}\sum\nolimits_{k \in [n]} L_{\mathbf{h}}\norm{a^i - a^i_\nu}_2 \leq L_{\mathbf{h}}\nu, \\
& \babs{\mathds{h}^i(a^i) - \mathds{h}^i(a^i_\nu)} \leq \bexpt{}{\norm{\mathbf{h}^i(a^i; X, Y) - \mathbf{h}^i(a^i_\nu; X, Y)}_2} \\ 
& \leq \bexpt{}{L_{\mathbf{h}} \norm{a^i - a^i_\nu}_2 } \leq L_{\mathbf{h}} \nu. 
\end{align*}
Then we have the following relations for the event under study: 
\begin{align*}
& \big\{\mathbb{A}_{-\varepsilon} \subseteq A_{\omega, n} \subseteq \mathbb{A}_{\varepsilon} \big\}^c \\
& \subseteq \big\{ \exists \bsb{a} \in \cls \mathbb{B}^N_R, i \in [N] \text{ s.t. } \abs{\mathbf{h}^i_{\omega, n}(a^i) - \mathds{h}^i(a^i)} > \varepsilon \big\} \\
& \subseteq \big\{ \exists \bsb{a} \in \cls \mathbb{B}^N_R, i \in [N] \text{ s.t. } \abs{\mathbf{h}^i_{\omega, n}(a^i_\nu) - \mathds{h}^i(a^i_\nu)} + 2L_{\mathbf{h}}\nu > \varepsilon \\
& \qquad \text{ and some } \bsb{a}_\nu \in A_\nu \text{ s.t. } \norm{\bsb{a} - \bsb{a}_\nu}_2 < \nu \big\} \\
& \subseteq \big\{ \exists \bsb{a}_\nu \in A_\nu, i \in [N] \text{ s.t. } \abs{\mathbf{h}^i_{\omega, n}(a^i_\nu) - \mathds{h}^i(a^i_\nu)} + 2L_{\mathbf{h}}\nu > \varepsilon \big\} \\
& = \big\{ \exists \bsb{a}_\nu \in A_\nu, i \in [N] \text{ s.t. } \abs{\mathbf{h}^i_{\omega, n}(a^i_\nu) - \mathds{h}^i(a^i_\nu)} > \varepsilon/2 \big\},
\end{align*}
where the second inclusion follows from the Lipschitz continuity derived above and the last inclusion is a result of applying the condition that $\nu = \varepsilon/(4L_{\mathbf{h}})$. 
By applying Lemma~\ref{le:set-approx-finite}, we arrive at the conclusion that 
\begin{align*}
\probb\big\{\mathbb{A}_{-\varepsilon} \subseteq A_{\omega, n} \subseteq \mathbb{A}_{\varepsilon} \big\}^c \leq 2N\mathfrak{n}_\nu  \exp\Big(- \frac{n\varepsilon^2}{12M_h^2}\Big).
\end{align*}

\subsubsection{Proof of Theorem~\ref{thm:as-res}}

Let $M \coloneqq (\sqrt{Nd}M_F) \vee M_h$. 
For each $n \in \nset{}{+}$, we have 
\begin{align*}
& \probb\{\mathscr{E}_n^c\} \leq \probb\big\{ \omega \mid \bnorm{\mathbf{F}_{\omega, n}(\bsb{a}) - \mathbb{F}(\bsb{a})}_2 < \varepsilon_n, \forall \bsb{a} \in \cls\mathbb{B}^N_R \big\}^c \\
& \qquad + \probb\big\{ \omega \mid \mathbb{A}_{-\varepsilon_n} \subseteq A_{\omega,n} \subseteq \mathbb{A}_{\varepsilon_n} \big\}^c \\
& \leq 2(Nd + N)\big(\frac{4LR\sqrt{Nd}}{\varepsilon_\diamond}n^{\alpha}\big)^{Nd}\exp\Big(- \frac{\epsilon_\diamond^2n^{1 - 2\alpha}}{12M^2}\Big) \\
& \overset{(a)}{=} C_1n^{\alpha Nd} \cdot \exp\big(-C_2 n^{1 - 2\alpha}\big) \overset{(b)}{\leq} \exp\big(-\frac{C_2}{2} n^{1 - 2\alpha}\big), 
\end{align*}
where in $(a)$, we set the constants $C_1$ and $C_2$ to be equal to the terms in the line above; 
there exists an index $n_\diamond$ such that $(b)$ holds for all $n > n_\diamond$. 
Then, it follows that 
\begin{align*}
\sum_{n \in \nset{}{+}} \probb\{\mathscr{E}_n^c\} \leq n_\diamond + \sum_{n > n_\diamond}\exp\big(-\frac{C_2}{2} n^{1 - 2\alpha}\big) < \infty.   
\end{align*}
By the Borel-Cantelli lemma \cite[Thm.~2.3.1]{durrett2019probability}, we have $\probb\big\{ \mathscr{E}_n^c \text{ i.o. } \big\} = 0$. 

\subsection{Proof for Solution Analysis}

\subsubsection{Proof of Lemma~\ref{le:cont-proj}}

Following the results of Theorem~\ref{thm:as-res}, for all $n > n_{\omega}$, it holds that $\mathbb{A}_{-\varepsilon_n} \subseteq A_{\omega, n} \subseteq \mathbb{A}_{\varepsilon_n}$, for $\varepsilon_n = \varepsilon_\diamond n^{-\alpha}$ with $0 < \alpha < 1/2$. 
We can observe that $\lim_{n \to \infty} \mathbb{A}_{-\varepsilon_n} = \mcup_{n \in \nset{}{+}} \mathbb{A}_{-\varepsilon_n} = \mathbb{A} \cap \{\bsb{a} \mid \mathds{h}(\bsb{a}) < 0\}$ and $\lim_{n \to \infty} \mathbb{A}_{\varepsilon_n} = \mcap_{n \in \nset{}{+}} \mathbb{A}_{\varepsilon_n} = \mathbb{A}$. 
Additionally, $A_{\omega, n}$ is a closed set for each $n \in \nset{}{+}$. 
Thus, $\lim_{n \to \infty} A_{\omega, n} = \mathbb{A}$. 

Regarding the second part of the lemma, we aim to prove the existence of the index $N_{0, \omega}$ and the open and bounded neighborhood $\Tilde{\mathcal{N}}$ of $\Tilde{\bsb{a}}_\sharp$. 
It suffices to prove that for an arbitrary sequence $(\Tilde{\bsb{a}}_n)_{n \in \nset{}{+}}$ converging to $\Tilde{\bsb{a}}_\sharp$, the sequence $(\proj_{A_{\omega,n}}[\Tilde{\bsb{a}}_n])_{n \in \nset{}{+}}$ converges to $\bsb{a}_\sharp$. 
This claim can be verified via \cite[Lemma~2.8.2]{facchinei2003finite}. 

\subsubsection{Proof of Lemma~\ref{le:unit-ind}}

We define the operator $\mathbf{F}_t(\bsb{a}) \coloneqq t \cdot \mathbb{F}(\bsb{a}) + (1-t)(\bsb{a} - \Tilde{\bsb{a}}_\star)$ with $t \in [0, 1]$. 
Then its normal map gives rise to a homotopy $H: \cls \mathbb{B}^N_R \times [0,1] \to \rset{Nd}{}$: 
\begin{align*}
& H(\Tilde{\bsb{a}}, t)= \mathbf{F}_t \circ \proj_{\mathbb{A}}[\Tilde{\bsb{a}}] + \Tilde{\bsb{a}}- \proj_{\mathbb{A}}[\Tilde{\bsb{a}}] \\
& = t \big( \mathbb{F} \circ \proj_{\mathbb{A}}[\Tilde{\bsb{a}}] + \Tilde{\bsb{a}} - \proj_{\mathbb{A}}[\Tilde{\bsb{a}}] \big) + (1-t)(\Tilde{\bsb{a}} - \Tilde{\bsb{a}}_\star) \\
& = t \mathbb{F}^{\text{nor}}(\Tilde{\bsb{a}}) + (1-t)(\Tilde{\bsb{a}} - \Tilde{\bsb{a}}_\star). 
\end{align*}
Let $\mathbb{B}_\varepsilon(\Tilde{\bsb{a}}_\star)$ denote the open ball in $\rset{Nd}{}$ centered at $\Tilde{\bsb{a}}_\star$ with radius $\varepsilon$. 
Assume ad absurdum that there exists a sequence of radius $(\varepsilon_k)_{k \in \nset{}{+}}$ such that $\varepsilon_k \searrow 0$ and for each $k \in \nset{}{+}$, there exists a solution pair $(\Tilde{\bsb{a}}_k, t_k) \in \partial\mathbb{B}_{\varepsilon_k}(\Tilde{\bsb{a}}_\star) \times [0,1]$ s.t. $H(\Tilde{\bsb{a}}_k, t_k) = 0$. 
This condition implies that $\bsb{a}_k \coloneqq \proj_{\mathbb{A}}[\Tilde{\bsb{a}}_k]$ is a solution of $\text{VI}(\mathbb{A}, \mathbf{F}_{t_k})$, which follows that for $\bsb{a}_\star$
\begin{align*}
& \langle {\bsb{a}}_\star - {\bsb{a}}_k, \mathbf{F}_{t_k}({\bsb{a}}_k)\rangle \geq 0 \iff\\
& t_k\langle {\bsb{a}}_\star - {\bsb{a}}_k, \mathbb{F}({\bsb{a}}_k)\rangle + (1 - t_k) \langle {\bsb{a}}_\star - {\bsb{a}}_k, \bsb{a}_k - \Tilde{\bsb{a}}_\star \rangle \geq 0\\
& \implies t_k\langle {\bsb{a}}_\star - {\bsb{a}}_k, \mathbb{F}({\bsb{a}}_k)\rangle - (1 - t_k)\norm{{\bsb{a}}_{k} - {\bsb{a}}_\star}_2^2 \geq 0,   
\end{align*}
where the last line follows from the observation $\langle {\bsb{a}}_\star - {\bsb{a}}_k, \bsb{a}_\star - \Tilde{\bsb{a}}_\star \rangle = \langle {\bsb{a}}_\star - {\bsb{a}}_k, \mathbb{F}(\bsb{a}_\star) \rangle \leq 0. $
Suppose $\bsb{a}_\star = \bsb{a}_k$, then $\Tilde{\bsb{a}}_k = \bsb{a}_k - \mathbf{F}_{t_k}(\bsb{a}_k) = \bsb{a}_\star - \mathbb{F}(\bsb{a}_\star) = \Tilde{\bsb{a}}_\star$, contradicting the setup $\Tilde{\bsb{a}}_k \in \partial \mathbb{B}_{\varepsilon_k}(\Tilde{\bsb{a}}_\star)$. 
Hence, $\bsb{a}_\star \neq \bsb{a}_k$ and $t_k \in (0, 1]$, which further implies that 
\begin{align*}
\langle {\bsb{a}}_\star - {\bsb{a}}_k, \mathbb{F}({\bsb{a}}_k)\rangle \geq \frac{1 - t_k}{t_k}\norm{{\bsb{a}}_{k} - {\bsb{a}}_\star}_2^2 \geq 0. 
\end{align*}
For another thing, $\mathbb{F}$ is $\lambda$-strongly monotone, i.e., $\langle \bsb{a}_k - \bsb{a}_\star, \mathbb{F}(\bsb{a}_k)\rangle \geq \lambda \norm{\bsb{a}_k - \bsb{a}_\star}^2_2 > 0$, leading to a contradiction. 
Hence, there exists an $\varepsilon > 0$ such that $\forall t \in [0, 1]$ and $\forall \Tilde{\bsb{a}} \in \partial \mathbb{B}_\varepsilon(\Tilde{\bsb{a}}_\star)$, $H(\Tilde{\bsb{a}}, t) \neq 0$. 
The homotopy invariance property of the degree (Theorem~\ref{thm:deg} (d3)) yields: 
\begin{align*}
& \ind(\mathbb{F}^{\text{nor}}, \Tilde{\bsb{a}}_\star) = \deg(\mathbb{F}^{\text{nor}}, \mathbb{B}_{\varepsilon}(\Tilde{\bsb{a}}_\star), 0) \\
& = \deg(H(\cdot, 1), \mathbb{B}_{\varepsilon}(\Tilde{\bsb{a}}_\star), 0) = \deg(H(\cdot, 0), \mathbb{B}_{\varepsilon}(\Tilde{\bsb{a}}_\star), 0) = 1.
\end{align*}

\subsubsection{Proof of Theorem~\ref{thm:sol-consis}}
Our focus is confined to the sample set $\Tilde{\Omega}$.  
As per the results in Lemma~\ref{le:cont-proj}, for any $\omega \in \Tilde{\Omega}$ and an arbitrary open and bounded neighborhood $\mathcal{N}$ of the unique solution $\bsb{a}_\star$, there exists an index $N_{0,\omega}$ sufficiently large and an open and bounded neighborhood $\Tilde{\mathcal{N}}$ of $\Tilde{\bsb{a}}_\star$ s.t. $\proj_{A_{\omega, n}}\big(\cls \Tilde{\mathcal{N}}\big) \subseteq \mathcal{N}$ for all $n > N_{0,\omega}$. 
Since $\mathbb{F}^{\text{nor}}$ is continuous and $\partial \Tilde{\mathcal{N}}$ is a compact set in $\rset{Nd}{}$, $\mathbb{F}^{\text{nor}}$ is a closed map and the image $\mathbb{F}^{\text{nor}}(\partial \Tilde{\mathcal{N}})$ is also closed. 
Moreover, $\Tilde{\bsb{a}}_\star \notin \partial \Tilde{\mathcal{N}}$, which implies that $\inf_{\Tilde{\bsb{a}} \in \partial \Tilde{\mathcal{N}}} \norm{\mathbb{F}^{\text{nor}}(\Tilde{\bsb{a}}) - \mathbb{F}^{\text{nor}}(\Tilde{\bsb{a}}_\star)}_2 = \inf_{\Tilde{\bsb{a}} \in \partial \Tilde{\mathcal{N}}} \norm{\mathbb{F}^{\text{nor}}(\Tilde{\bsb{a}})}_2 = \ubar{\delta} > 0$. 
Subsequently, we aim to demonstrate that the discrepancy between $\mathbb{F}^{\text{nor}}$ and $\mathbf{F}^{\text{nor}}_{\omega, n}$ within the vicinity of $\Tilde{\bsb{a}}_\star$ is mitigable for $n$ sufficiently large. 
For any $\Tilde{\bsb{a}} \in \cls \Tilde{\mathcal{N}}$ and $n > n_\omega$ with the index $n_{\omega}$ defined in Theorem~\ref{thm:as-res}, we have
\begin{align*}
& \norm{\mathbb{F}^{\text{nor}}(\Tilde{\bsb{a}}) - \mathbf{F}^{\text{nor}}_{\omega, n}(\Tilde{\bsb{a}})}_2 = \norm{\mathbf{F}_{\omega, n}\circ \proj_{A_{\omega,n}}(\Tilde{\bsb{a}})  - \mathbf{F}_{\omega, n} \circ \proj_{\mathbb{A}}(\Tilde{\bsb{a}}) \\
& \quad  + \mathbf{F}_{\omega, n} \circ P_{\mathbb{A}}(\Tilde{\bsb{a}}) - \mathbb{F} \circ \proj_{\mathbb{A}}(\Tilde{\bsb{a}}) - \proj_{A_{\omega, n}}(\Tilde{\bsb{a}}) + \proj_{\mathbb{A}}(\Tilde{\bsb{a}})}_2 \\
& \overset{(a)}{\leq} (L_F + 1)\norm{\proj_{A_{\omega, n}}(\Tilde{\bsb{a}}) - \proj_{\mathbb{A}}(\Tilde{\bsb{a}})}_2 + \norm{(\mathbf{F}_{\omega, n} - \mathbb{F}) \circ \proj_{\mathbb{A}}(\Tilde{\bsb{a}})}_2 \\
& \overset{(b)}{\leq} (L_F + 1)\norm{\proj_{A_{\omega, n}}(\Tilde{\bsb{a}}) - \proj_{\mathbb{A}}(\Tilde{\bsb{a}})}_2 + \varepsilon_\diamond n^{-\alpha}, 
\end{align*}
where in $(a)$, we apply the Lipschitz continuity derived in the proof of Proposition~\ref{propst:opt-approx-cmpct}; 
$(b)$ follows from the implication of Theorem~\ref{thm:as-res}. 
Again, we invoke Lemma~\ref{le:cont-proj}, obtaining that for any $\Tilde{\bsb{a}} \in \cls \Tilde{\mathcal{N}}$, $\lim_{n \to \infty} \norm{\proj_{A_{\omega, n}}(\Tilde{\bsb{a}}) - \proj_{\mathbb{A}}(\Tilde{\bsb{a}})}_2 = 0$. 
By the fact that $\cls \Tilde{\mathcal{N}}$ is a compact set, together with the previous observations, there exists an index $N_{\omega} > N_{0, \omega}$, such that $\norm{\mathbb{F}^{\text{nor}}(\Tilde{\bsb{a}}) - \mathbf{F}^{\text{nor}}_{\omega, n}(\Tilde{\bsb{a}})}_2 < \ubar{\delta}$ for all $\Tilde{\bsb{a}} \in \cls\Tilde{\mathcal{N}}$. 
It further implies that $\norm{\mathbf{F}^{\text{nor}}_{\omega, n}(\Tilde{\bsb{a}})}_2 \geq \norm{\mathbb{F}^{\text{nor}}(\Tilde{\bsb{a}})}_2 - \norm{\mathbb{F}^{\text{nor}}(\Tilde{\bsb{a}}) - \mathbf{F}^{\text{nor}}_{\omega, n}(\Tilde{\bsb{a}})}_2 \geq \ubar{\delta} - \norm{\mathbb{F}^{\text{nor}}(\Tilde{\bsb{a}}) - \mathbf{F}^{\text{nor}}_{\omega, n}(\Tilde{\bsb{a}})}_2 > 0$ for all $\Tilde{\bsb{a}} \in \partial \Tilde{\mathcal{N}}$. 

Construct the homotopy $H(\Tilde{\bsb{a}}, t) = t\mathbf{F}^{\text{nor}}_{\omega, n}(\Tilde{\bsb{a}}) + (1 - t)\mathbb{F}^{\text{nor}}(\Tilde{\bsb{a}})$ and fix $y = 0$. 
In light of the above reasoning, $y \not \in H(\partial \Tilde{\mathcal{N}}, t)$ for all $t \in [0, 1]$. 
By the homotopy invariance property of the degree (Theorem~\ref{thm:deg} (d3)), we have
\begin{align*}
\deg(\mathbf{F}^{\text{nor}}_{\omega, n}, \Tilde{\mathcal{N}}, 0) 
= \deg(\mathbb{F}^{\text{nor}}, \Tilde{\mathcal{N}}, 0) = \ind(\mathbb{F}^{\text{nor}}, \Tilde{\bsb{a}}_\star) = 1.
\end{align*}
As per Theorem~\ref{thm:deg}, the preimage $\mathbf{F}_{\omega, n}^{-1}(0)$ intersects $\Tilde{\mathcal{N}}$, i.e., $\Tilde{\bsb{a}}_{\omega, n} \in \Tilde{\mathcal{N}}$ and $\bsb{a}_{\omega, n} \in \mathcal{N}$. 
For an arbitrary $\varepsilon > 0$, choosing $\mathcal{N} = \mathbb{B}_{\varepsilon}(\bsb{a}_\star)$ yields that $\lim_{n \to \infty} \norm{\bsb{a}_{\omega,n} - \bsb{a}_\star}_2 = 0$.

\ifarxiv
\subsection{Numerical Experiments}

In this section, we focus on solving a multi-account portfolio optimization problem, a variant of \cite{yang2013multi}. 
The game $\mathcal{G}$ involves $N$ investment accounts/players and $m$ financial assets.   
When portfolio managers make investment decisions, they are unable to directly access the actual return rates $Y: \Omega \to \mathcal{Y}$ associated with assets. 
Instead, they can monitor global factors $X: \Omega \to \mathcal{X}$, such as market demands, unemployment rates, inflation rates, and demographic trends, to determine the most suitable asset combinations. 

Here, $\mathcal{X} \subseteq \rset{n_x}{}$ and $\mathcal{Y} \subseteq \rset{n_y}{}$. 
For each $i \in [N], x \in \mathcal{X}$, the decision function is represented as $\hat{z}^i(x) = [\hat{b}^i(x); \hat{s}^i(x)] \in \rset{2n_y}{+}$, where $\hat{b}^i(x)$ and $\hat{s}^i(x)$ correspond to the purchase and sale amount of the assets, when the observation $X = x$ is made. 
Furthermore, we also consider the market impact cost function $\psi:\rset{2n_y}{+} \to \rset{2n_y}{+}$ that induces the costs associated with the rebalancing of the current investment combinations to new ones. 
Let $b = \sum_{i \in [N]} b^i \in \rset{n_y}{+}$ and $s = \sum_{i \in [N]} s^i \in \rset{n_y}{+}$ denote the amounts of collective buying and selling. 
Then, the market impact cost is defined as 
$\psi([b; s]) = \Big[\begin{smallmatrix}\Omega^+ & 0 \\ 0 & \Omega^- \end{smallmatrix}\Big] \Big[\begin{smallmatrix} b \\ s \end{smallmatrix}\Big]$, where $\Omega^+ \in \rset{n_y \times n_y}{}$ 
and $\Omega^- \in \rset{n_y \times n_y}{}$ are two positive definite diagonal matrices. 
The cost function of player $i$ is given by 
\begin{align*}
c^i(z^i; y, z^{-i}) \coloneqq - \begin{bmatrix} y & -y \end{bmatrix}z^i + (z^i)^T\psi \big(\sum\nolimits_{j \in [N]}z^j \big). 
\end{align*}
To enforce the constraints that each entry of $z^i$ should be non-negative and the funding budget $[\bsb{1}^T, -\bsb{1}^T]z^i \leq w^i$ for the budget upper bound $w^i > 0$, we further add the penalty term to the cost function 
\begin{align*}
p^i(z^i) \coloneqq \lambda'\Big( \sum\nolimits_{l \in [2n_y]}\big[-[z^i]_l\big]_+^2 + \big[[\bsb{1}^T, -\bsb{1}^T]z^i - w^i\big]_+^2\Big), 
\end{align*}
where $\lambda' > 0$ is some weighting coefficient and $[x]_+$ denotes a smoothed version of $\max(0, x)$. 
The objective function is given by $J^i(z^i; y, z^{-i}) \coloneqq c^i(z^i; y, z^{-i}) + p^i(z^i)$. 
When optimizing the augmented objective function, each player $i$ also considers the following CVaR expectation constraint \cite{rockafellar2002conditional}: 
\begin{align*}
&  \min_{\zeta^i} \Big\{h^i(z^i, \zeta^i; y) \Big\} \leq 0, \text{ with } h^i(z^i, \zeta^i; y) \coloneqq \\
&  \qquad \zeta^i + \frac{1}{1-\alpha^i}\Bexpt{X, Y}{\big[-[Y;-Y]^T z^i(X) - \zeta^i\big]_+} - \gamma^i, 
\end{align*}
which mandates that the cost remains below $\gamma^i$ with a probability greater than $\alpha^i$. 
With the inclusion of extra regularization terms to promote solution consistency, the local optimization problem of player $i$ can be formulated as follows: 
\begin{align*}
& a^i_\star = \underset{a^i \in \mathbb{A}^i, \zeta^i \in \rset{}{}}{\argmin} \bexpt{}{\mathbf{J}^i\big(a^i; X, Y, a^{-i}\big)} + \frac{\lambda}{2}\big(\norm{a^i}^2_2 + \abs{\zeta^i}^2\big) \\
& \mathbb{A}^i \coloneqq \{a^i \in \cls\mathbb{B}_R \mid \mathds{h}^i(a^i, \zeta^i) = 
 \expt{}{\mathbf{h}^i(a^i, \zeta^i; X, Y)} \leq 0\}. 
\end{align*}
According to \cite[Sec.~4.2]{lan2012validation}, $\zeta^i$ can be confined to a bounded interval, preserving the compactness of the feasible set. 

In our experiments, we consider $N = 3$ investment accounts, $n_y = 4$ assets, and $n_x = 5$ observable features. 
The actual return rate $Y$ is modeled as a quadratic function of $X$, accompanied by additional random factors with non-zero means. 
The set of basis functions $\Phi^d$ consists of monomials in $x_1, \ldots, x_5$ with degrees no greater than $2$ and $\abs{\Phi^d} = 21$.
For the market impact cost, we specify that $\Omega^- = 4\Omega^+$, as a result of which the cost incurred from selling assets surpasses that arising from purchasing them. 
We let $\lambda' = 50$ and $\lambda = 10^{-4}$. 
For the three players involved, we assign $\alpha^i$ values as $(0.95, 0.8, 0.8)$, $\gamma^i$ values as $(-1, 5, 5)$, and $w^i$ values as $(10, 40, 2)$.
This assignment characterizes player $1$ as a conservative investor with moderate funds, player $2$ as an aggressive investor with ample funds, and player $3$ as an aggressive investor with limited funds. 
In the offline learning setup, we consider training datasets of size $n \in \{10, 20, 100, 200, 10^3, 2\times 10^3\}$ and evaluate the performance of the decision functions obtained through a test dataset of size $n_{\text{test}} = 10^4$. 
The finite-sample approximation problem \eqref{eq:prob-fs} is solved via the extra-gradient algorithm with a constant step size $3 \times 10^{-4}$ and a total number of $10^5$ iterations. 
In the pursuit of a reference solution $\bsb{a}_\star$, the extra-gradient algorithm working in an online-learning manner 
\cite{iusem2017extragradient}, is implemented. 
Achieving convergence of the algorithm proposed in \cite{iusem2017extragradient} requires superlinearly growing sample counts over the course of iterations, and our implementation leverages over $2\times 10^8$ samples to fulfill the superlinearly growing requirement.  

In the upper panel of Fig.~\ref{fig:payoffs}, we visualize the average payoff function values in the test dataset, that is, $-c^i$, and the augmented payoff function values that incorporate penalties for violation of constraints $\bar{p}^i$ and regularization terms. 
The values of the augmented payoff function are depicted by blue patches, while the orange patches illustrate $p^i$ and the associated regularization terms. 
The sum of these components results in $-c^i$. 
The green line illustrates the evolution of the relative distance $\norm{\bsb{a}_{\omega, n} - \bsb{a}_\star}_2/\norm{\bsb{a}_\star}_2$ with increasing values of $n$. 
As $n$ increases, the relative distance shows a steady decrease, and simultaneously, penalties associated with constraint violations and regularization exhibit a decrease. 
Player $2$ achieves higher investment income compared to the conservative investor, player $1$, and player $3$, who has a limited investment budget.
Table~\ref{tab:constr} provides the average values of the constraint function obtained from the testing set. 
Notably, all finite-sample approximations of the constraints are fulfilled except for player $1$ when $n = 10$.
This deviation is expected, given the conservative investment strategy of player $1$ ($\alpha^1 = 0.95$, $\gamma^1 = -1$) and the limited number of training samples $n = 10$. 
To demonstrate the sample efficiency of the offline learning approach compared to the online ones, the lower panel of Fig.~\ref{fig:payoffs} illustrates the relative distances between the coefficients $\bsb{a}_{\omega, n}$ obtained via \cite[Algorithm~1]{iusem2017extragradient} and the solution $\bsb{a}_\star$ for varying step sizes.
We denote the stepsize in the implementation by $\tau$, and increasing $\tau$ further will cause the sequence to become divergent.
In the online learning paradigm, despite processing more than $5 \times 10^5$ samples, the coefficients $\bsb{a}_{\omega, n}$ continue to show a relative distance greater than 0.5 from the solution $\bsb{a}_\star$.

\begin{figure}
    \centering
    \includegraphics[width=0.45\textwidth]{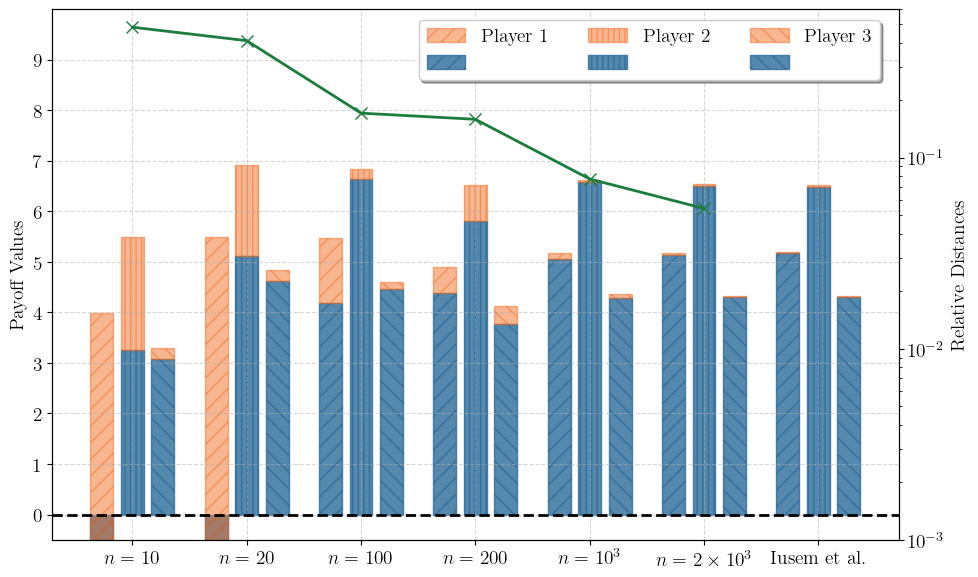}
    \includegraphics[width=0.45\textwidth]{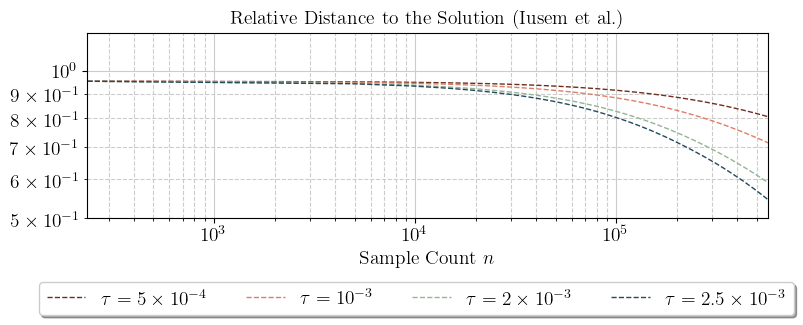}
    \caption{Analyzing the Influence of Training Sample Size on Payoffs and Relative Distances to the Solution}
    \label{fig:payoffs}
\end{figure}
\begin{table}
    \centering
    \begin{tabular}{cccc}
        \toprule
        & Player 1 & Player 2 & Player 3 \\
        \midrule
        $n=10$ & 0.4051 & -4.7171 & -4.7652 \\
        $n=20$ & -0.2311 & -4.8321 & -4.8149 \\
        $n=100$ & -0.4348 & -4.8186 & -4.8224 \\
        $n=200$ & -0.6422 & -4.7833 & -4.8128 \\
        $n=10^3$ & -0.2970 & -4.7976 & -4.8186 \\
        $n=2 \times 10^3$ & -0.2027 & -4.7801 & -4.8166 \\
        Iusem et al. \cite{iusem2017extragradient} & -0.3295 & -4.7976 & -4.8188 \\
        \bottomrule
    \end{tabular}
    \caption{Results of Constraint Function Evaluations}
    \label{tab:constr}
    \vspace*{-\baselineskip}
\end{table}
\fi

\end{document}